\documentclass[12pt]{article}
\usepackage{graphicx}
\usepackage{amsmath}
\usepackage{amsfonts}
\usepackage{amssymb}
\newtheorem{theorem}{Theorem}

\newtheorem{corollary}[theorem]{Corollary}

\newtheorem{lemma}[theorem]{Lemma}

\newtheorem{proposition}[theorem]{Proposition}

\newenvironment{proof}[1][Proof]{\textbf{#1.} }{\ \rule{0.5em}{0.5em}}
\def\text{\hbox} 

\def\a{\alpha}
\def\b{\beta}
\def\bb{\gamma}

\def\p{\pi}

\def\vf{\psi}

\def\f{\varphi}

\def\S{{\bf S}}

\def\A{{\cal A}}
\def\B{{\cal B}}
\def\K{{\cal K}}
\def\A'{{\cal A}'}

\date{}

\begin{document}

\title{Torus knots and Dunwoody manifolds}

\author{ H\"{u}seyin Aydin \and Inci G\"{u}ltekyn \and Michele Mulazzani}

\maketitle

%\centerline{(preliminary version)}

\begin{abstract}
We obtain an explicit representation, as Dunwoody manifolds, of
all cyclic branched coverings of torus knots of type $(p,mp\pm
1)$, with $p>1$ and $m>0$.\\
\\ {{\it Mathematics Subject
Classification 2000:} Primary 57M05, 20F38; Secondary 57M12, 57M25.}

\end{abstract}

%\newpage

\section{Introduction} \label{intro}

Many interesting examples of cyclic branched coverings of knots in
$\S^3$ admitting cyclic presentations for the fundamental groups
have recently been found (see \cite{CHK,Du,HKM,Ki,KKV1,MR,VK}). In
order to investigate these relations, M. J. Dunwoody introduced in
\cite{Du} a class of 3-manifolds, depending on six integer
parameters, with cyclically presented fundamental groups. It has
been shown in \cite{GM} that all these manifolds turn out to be
strongly-cyclic branched coverings of $(1,1)$-knots in lens spaces
(possibly $\S^3$). Moreover, the explicit Dunwoody representation
for all cyclic branched coverings of 2-bridge knots has been
obtained.

In this paper we give a similar result for a wide class of torus
knots, which are, together with 2-bridge knots, the most important
examples of $(1,1)$-knots in $\S^3$. The Dunwoody parameters are
obtained for all cyclic branched coverings of torus knots of type
$(p,mp\pm 1)$, with $p>1$ and $m>0$, thus including all torus
knots with bridge number $\le 4$. These manifolds have been
considered, from a different point of view, in \cite{CM1} and
\cite{CHK}.

We refer to \cite{Ka,Ro} for details on knot theory and cyclic
branched coverings of knots and to \cite{Jo} for details on cyclic
presentation of groups.

\section{(1,1)-knots and Dunwoody manifolds}

A knot $K$ in a 3-manifold $N^3$ is called a $(1,1)$-{\it knot\/}
if there exists a Heegaard splitting of genus one
$$(N^3,K)=(T,A)\cup_{\f}(T',A'),$$ where $T$ and $T'$ are solid
tori, $A\subset T$ and $A'\subset T'$ are properly embedded
trivial arcs, and $\f:(\partial T',\partial A')\to(\partial
T,\partial A)$ is an attaching homeomorphism (see Figure \ref{Fig.
1}). Obviously, $N^3$ turns out to be a lens space $L(p,q)$
(including $\S^3=L(1,0)$).

%\bigskip

\begin{figure}[ht]
\begin{center}
\includegraphics*[totalheight=3cm]{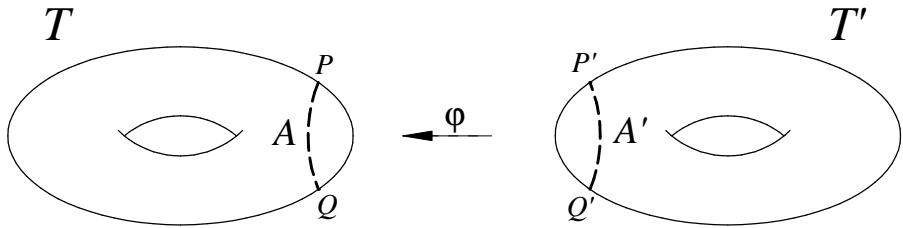}
\end{center}
\caption{A $(1,1)$-decomposition.} \label{Fig. 1}
\end{figure}

It is well known that the family of $(1,1)$-knots contains all
torus knots and all two-bridge knots in $\S^3$. Several
topological properties of $(1,1)$-knots have recently been
investigated (see references in \cite{CM2}).

An algebraic representation of $(1,1)$-knots has been developed in
\cite{CM1} and \cite{CM2}, where it is shown that there is a
natural surjective map $$\vf\in PMCG_2(\partial T)\mapsto
K_{\vf}\in \K_{1,1}$$ from the pure mapping class group of the
twice punctured torus $PMCG_2(\partial T)$ to the class $\K_{1,1}$
of all $(1,1)$-knots.

\begin{figure}[ht]
\begin{center}
\includegraphics*[totalheight=3.5cm]{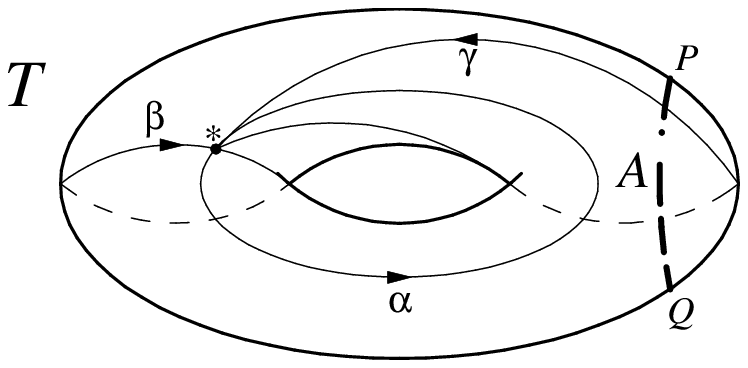}
\end{center}
\caption{} \label{Fig. 0}
\end{figure}

The fundamental group of the exterior of a $(1,1)$-knot $K_{\vf}$
can be explicitly obtained using its representation $\vf$. Let
$\a,\b,\bb\subset\partial T$ be the loops depicted in Figure
\ref{Fig. 0}. They represent a set of free generators for
$\pi_1(\partial T -\partial A,*)$, while $\a$ and $\bb$ freely
generate $\pi_1(T-A,*)$. A straightforward application of
Seifert-Van Kampen theorem gives:

\begin{lemma} \textup{\cite{CM1}}
\label{fundamental1} The fundamental group of the exterior of a
$(1,1)$-knot $K_{\vf}\subset L(p,q)$ admits the presentation
$$\p_1(L(p,q)-K_{\vf},*)=\langle\,\a,\bb\,|\,r(\a,\bb)\,\rangle,$$
where $r(\a,\bb)$ is the homotopy class of $i(\vf(\b))$, being
$i:(\partial T,\partial A)\to(T,A)$ the inclusion map.
\end{lemma}

As a consequence, all $(1,1)$-knots in $\S^3$ are prime \cite{No},
and can therefore be classified, up to mirror image, by their
fundamental groups (see \cite[p. 76]{Ka}). Many results concerning
the connections between cyclically presented groups and
strongly-cyclic branched coverings of $(1,1)$-knots (see
definition in \cite{CM1}) have been obtained. It has been proved
in \cite{Mu} that every $n$-fold strongly-cyclic branched covering
of a $(1,1)$-knot admits a Heegaard diagram of genus $n$ which
encodes a cyclic presentation for the fundamental group. This
result has been improved in \cite{CM1}, where a constructive
algorithm which explicitly gives the cyclic presentations is
obtained.

\begin{figure}
 \begin{center}
 \includegraphics*[totalheight=6.5cm]{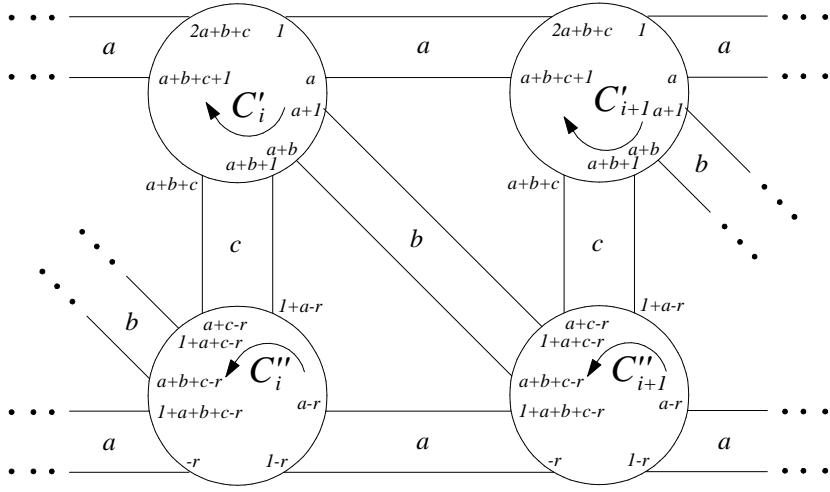}
 \end{center}
 \caption{A Dunwoody diagram.}

 \label{Fig. 13}

\end{figure}

\medskip

The family of Dunwoody manifolds has been introduced in \cite{Du}
by a class of trivalent regular planar graphs (called {\it
Dunwoody diagrams\/}) with cyclic symmetry, depending on six
integers $a,b,c,n,r,s$, such that $n>0$, $a,b,c\ge 0$ and
$a+b+c>0$. For certain values of the parameters, called {\it
admissible\/}, the Dunwoody diagrams $D(a,b,c,n,r,s)$ turn out to
be Heegaard diagrams, so defining a wide class of closed,
orientable 3-manifolds $M(a,b,c,n,r,s)$ with cyclically presented
fundamental groups.

More precisely, an admissible Dunwoody diagram $D(a,b,c,n,r,s)$ is
an open Heegaard diagram of genus $n$ (see Figure \ref{Fig. 13}),
which contains $n$ upper cycles $C'_1,\ldots,C'_n$, and $n$ lower
cycles $C''_1,\ldots,C''_n$, each having $d=2a+b+c$ vertices. For
every $i=1,\ldots,n$, the cycle $C'_i$ (resp. $C''_i$) is
connected to the cycle $C'_{i+1}$ (resp. $C''_{i+1}$) by $a$
parallel arcs, to the cycle $C''_{i}$ by $c$ parallel arcs and to
the cycle $C''_{i+1}$ by $b$ parallel arcs. We denote by $\B$ the
set of all these arcs. The cycle $C'_i$ is glued to the cycle
$C''_{i-s}$ (subscripts mod $n$) so that equally labelled vertices
are identified together. Observe that the parameters $r$ and $s$
can be considered mod $d$ and $n$ respectively. Since the
identification rule and the diagram are invariant with respect to
a cyclic action of order $n$, the Dunwoody manifolds admit a
cyclic symmetry of order $n$. Obviously, the Dunwoody manifold
$M(a,b,c,1,r,0)$ is homeomorphic to a lens space (possibly
$\S^3$), since it admits a Heegaard splitting of genus one.

A characterization of all Dunwoody manifolds as strongly-cyclic
branched coverings of $(1,1)$-knots is given by the following:

\begin{proposition} \label{strongly-cyclic} \textup{\cite{GM}}
The Dunwoody manifold $M(a,b,c,n,r,s)$ is the $n$-fold
strongly-cyclic covering of the lens space $M(a,b,c,1,r,0)$
(possibly $\S^3$), branched over a $(1,1)$-knot $K(a,b,c,r)$ only
depending on the integers $a,b,c,r$.
\end{proposition}

The converse of this result is also true, as proved in \cite{CM3}.
So the class of the Dunwoody manifolds coincides with the class of
the strongly-cyclic branched coverings of $(1,1)$-knots.

\begin{figure}
 \begin{center}
 \includegraphics*[totalheight=8.5cm]{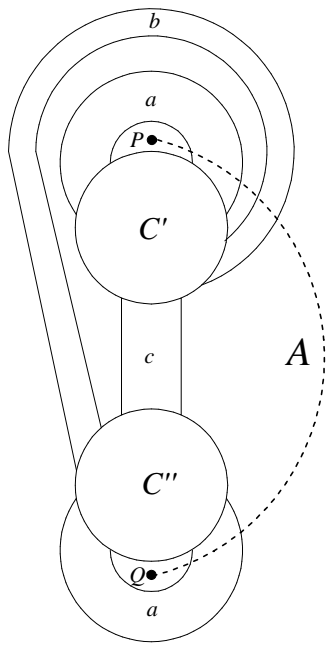}
 \end{center}
 \caption{$D(a,b,c,1,r,0)$}

 \label{Fig. 4}

\end{figure}

\begin{figure}
 \begin{center}
 \includegraphics*[totalheight=8.2cm]{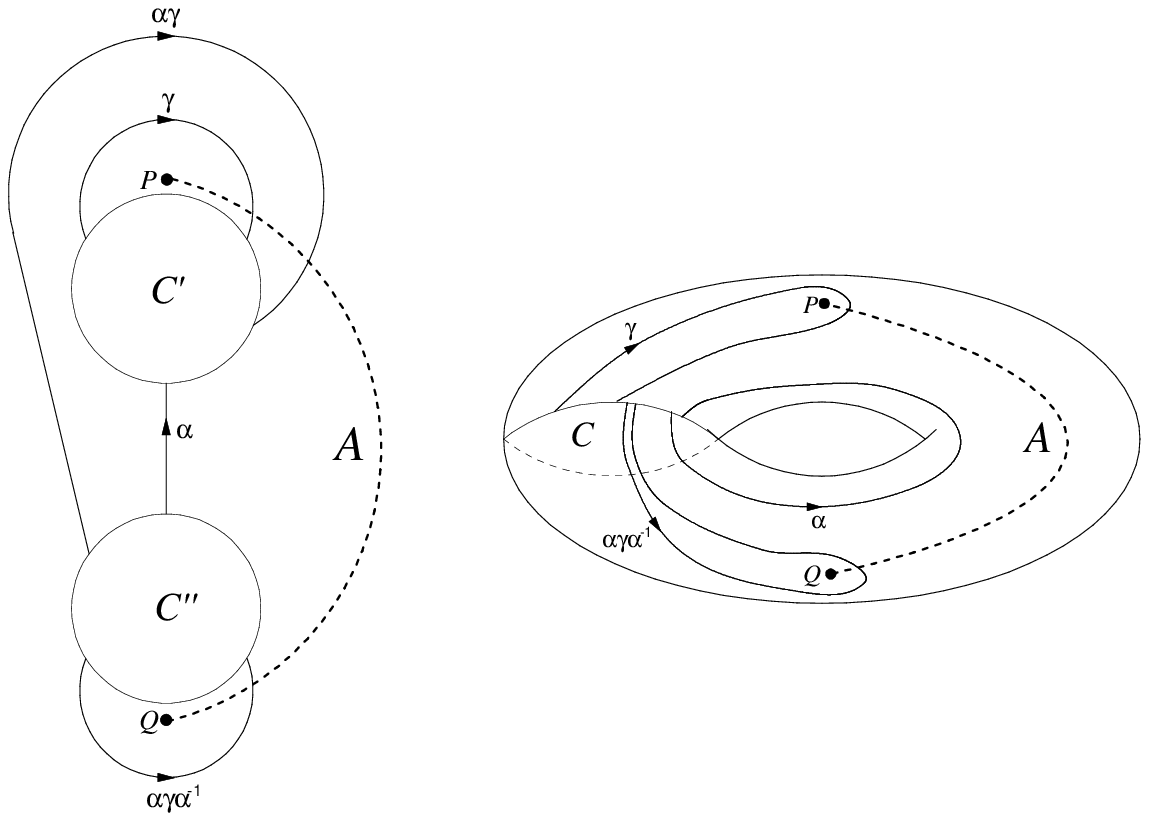}
 \end{center}
 \caption{}

 \label{Fig. 4'}

\end{figure}

The $(1,1)$-knots $K(a,b,c,r)$ occurring in Proposition
\ref{strongly-cyclic} admit a natural $(1,1)$-decomposition
$(T,A)\cup_{\f}(T',A')$ depicted in Figure \ref{Fig. 4}, where the
arcs of $\B$ constitute the curve $\f(\b')=\vf(\b)$, $\vf$ being
the element of $PMCG_2(\partial T)$ correspondent to $\f$. The
fundamental group of the exterior of $K(a,b,c,r)$ can be directly
read in the Dunwoody diagram of $D(a,b,c,1,r,0)$. The relation
$r(\a,\bb)$ of the presentation of Lemma \ref{fundamental1} is
obtained by walking along the arcs of $\B$, following a fixed
orientation: associate to each arc a word in $\a$ and $\bb$
representing its homotopy class in the fundamental group of $T-A$
(see Figure \ref{Fig. 4'}), where a properly embedded disk with
boundary $C$ is considered squeezed to the base point $*$ .

%\newpage

\section{Main results}

An interesting problem is to find the Dunwoody parameters of the
cyclic branched coverings of important classes of $(1,1)$-knots,
in particular when the knot lies in $\S^3$. This type of result
has been obtained in \cite[Theorem 8]{GM} for all 2-bridge knots.

Now we obtain a similar result for the torus knots ${\bf
t}(p,mp\pm 1)$, with $m>0$ and $p>1$.

\begin{figure}
 \begin{center}
 \includegraphics*[totalheight=12cm]{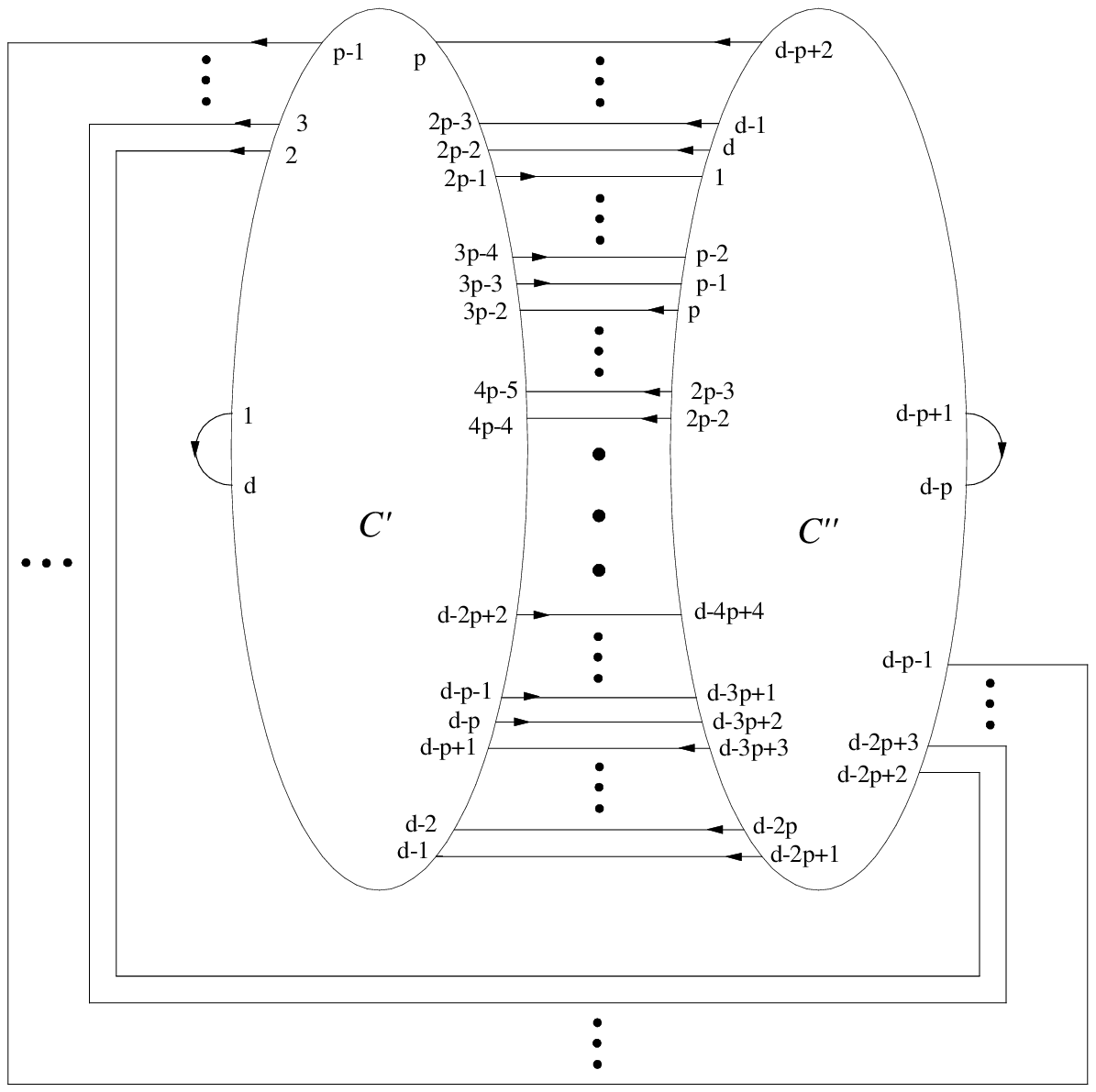}
 \end{center}
 \caption{$D(1,p-2,2mp-2m-p+1,1,p,0)$.}
 \label{Fig. 11}
\bigskip\bigskip
\end{figure}

\begin{figure}
 \begin{center}
 \includegraphics*[totalheight=17cm]{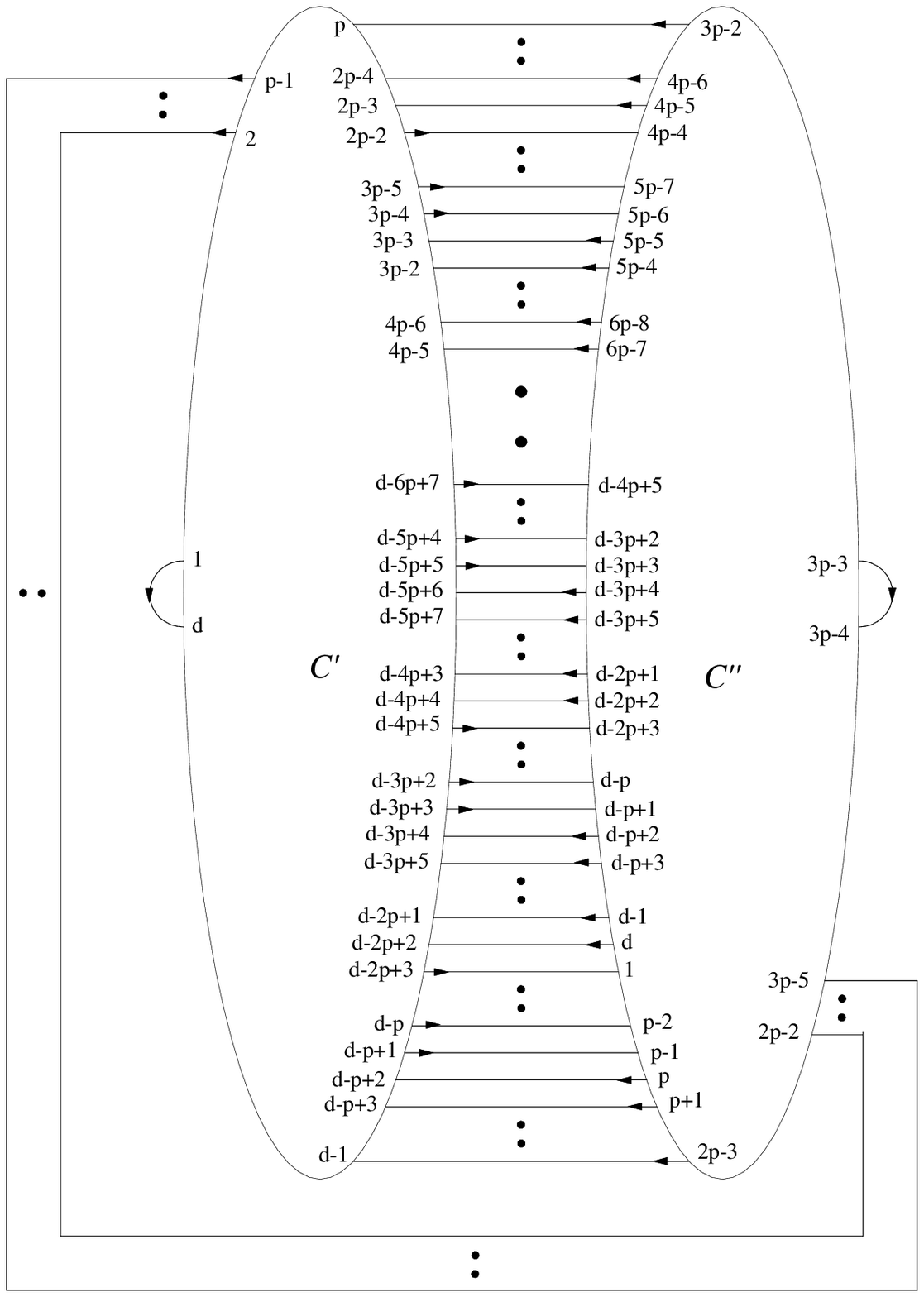}
 \end{center}
 \caption{$D(1,p-2,2mp-2m-p-1,1,-3p+4,0)$.}
 \label{Fig. 12}
\bigskip\bigskip
\end{figure}

\vbox{
\begin{proposition} \label{turchia+-}
\begin{enumerate}
\item $K(1,p-2,2mp-2m-p+1,p)$ is the torus knot ${\bf t}(p,mp+1)$,
for all $m>0$ and $p>1$;
\item $K(1,p-2,2mp-2m-p-1,-3p+4)$ is the
torus knot ${\bf t}(p,mp-1)$, for all $m>1$ and $p>1$.
\end{enumerate}
\end{proposition}
}
\begin{proof} {\it 1.} Let $K=K(1,p-2,2mp-2m-p+1,p)$.
Figure \ref{Fig. 11} depicts the Dunwoody diagram
$D=D(1,p-2,2mp-2m-p+1,1,p,0)$. The number of vertices of each
cycle is $d=2m(p-1)+1$. Starting from the vertex of $C''$ with
label $d$, we pass along all the arcs of $\B$ with the described
orientation, obtaining the word
$w=\a^{m(p-1)}\a\bb^{-1}\a^{-1}(\a^{-(m-1)}\bb^{-1}\a^{-1})^{p-2}\a^{-(m-1)}\bb^{-1}$.
Therefore, since the conditions of \cite[Corollary 4]{GM} are
satisfied, $D$ is admissible. Moreover, the fundamental group of
$M(1,p-2,2mp-2m-p+1,1,p,0)$ is $\langle \a,\bb\mid w,\bb\rangle$,
which is trivial. So $M(1,p-2,(p-1)(2m-1),1,p,0)\cong \S^3$.
Moreover, $\pi_1(\S^3-K)=\langle \a,\bb\mid w\rangle$. Since
$w=\a^{m(p-1)+1}(\bb^{-1}\a^{-m})^{p-1}\bb^{-1}=\a^{-m}\a^{mp+1}(\bb^{-1}\a^{-m})^p\a^m$,
we have $\pi_1(\S^3-K)\cong \langle \a,\bb\mid w'\rangle$, with
$w'=\a^{mp+1}(\bb^{-1}\a^{-m})^p$. Obviously this group is
isomorphic to the group $\langle x,y\mid x^{mp+1}y^{-p}\rangle$,
which is the group of the torus knot ${\bf t}(p,mp+1)$. As a
consequence $K$ is precisely ${\bf t}(p,mp+1)$.

{\it 2.} The proof is analogous to the previous one. Let
$K=K(1,p-2,2mp-2m-p-1,-3p+4)$. Figure \ref{Fig. 12} depicts the
Dunwoody diagram $D=D(1,p-2,2mp-2m-p-1,1,-3p+4,0)$. The number of
vertices of each cycle is $d=2m(p-1)-1$ (note that the labelling
is considered mod $d$). Starting from the vertex of $C''$ with
label $d$, we pass along all the arcs of $\B$, obtaining the word
$w=\a^{m-1}(\a^{m})^{p-3}\a^{m-1}\a\bb^{-1}\a^{-1}\a^{-(m-1)}(\bb^{-1}\a^{-1}\a^{-(m-1)})^{p-2}\bb^{-1}$.
So $D$ is admissible. Moreover, $M(1,p-2,2mp-2m-p-1,1,-3p+4,0)$ is
homeomorphic to $\S^3$, since its fundamental group is trivial.
Since
$w=\a^{m(p-1)-1}(\bb^{-1}\a^{-m})^{p-1}\bb^{-1}=\a^{mp-1}(\a^{-m}\bb^{-1})^p$,
the group $\pi_1(\S^3-K)\cong \langle \a,\bb\mid w\rangle$ is
isomorphic to the group $\langle x,y\mid x^{mp-1}y^{-p}\rangle$.
Therefore, $K$ is the torus knot ${\bf t}(p,mp-1)$.
\end{proof}

\vbox{
\begin{corollary} \label{turchia-n}
\begin{enumerate}
\item For all $m>0$ and $p>1$, the $n$-fold cyclic branched
covering of the torus knot ${\bf t}(p,mp+1)$ is the Dunwoody
manifold $M(1,p-2,2mp-2m-p+1,n,p,p)$; \item for all $m>1$ and
$p>1$, the $n$-fold cyclic branched covering of the torus knot
${\bf t}(p,mp-1)$ is the Dunwoody manifold
\hbox{$M(1,p-2,2mp-2m-p-1,n,-3p+4,-p)$.}
\end{enumerate}
\end{corollary}
}
\begin{proof}
From Proposition \ref{strongly-cyclic}, we only need to find the
sixth parameter $s$ of the 6-tuples. The relation
$q_{\sigma}+sp_{\sigma}\equiv 0$ mod $n$ must be satisfied for all
$n$, where $p_{\sigma}$ is the number of arcs of $\B$ oriented
from $C'$ to $C''$ minus the number of the ones oriented from the
$C''$ to $C'$ and $q_{\sigma}$ is the number of arcs of $\B$
oriented from the right to the left minus the number of the ones
oriented from the left to the right in the Dunwoody diagram (see
\cite[p. 385]{GM}). In other words, $-p_{\sigma}$ is the total
exponent of $\a$ in $w$ and $-q_{\sigma}$ is the total exponent of
$\bb$ in $w$. In the first case $q_{\sigma}=p$ and
$p_{\sigma}=-1$; therefore $s=p$. In the second case
$q_{\sigma}=p$ and $p_{\sigma}=1$; therefore $s=-p$.
\end{proof}

\medskip

An extension of the previous results to the whole class of torus
knots by using the same proof technique seems to be very
complicated, although some partial results have been obtained. A
possible alternative method could be to work by induction on the
number of steps of the Euclidean division algorithm needed to
obtain the greatest common division of the parameters of the torus
knot. Following this approach, the results of this article could
represent the first inductive step.

\bigskip\bigskip

\noindent {\bf Acknowledgements.} Work partially developed during
a visit by the third named author to the Department of Mathematics
of Atat\"{u}rk University, Erzurum (Turkey), supported by
TOKTEN/UNISTAR funds, and performed under the auspices of the
G.N.S.A.G.A. of I.N.d.A.M. (Italy) and the University of Bologna
funds for selected research topics.

%\newpage

%\newpage
\bigskip\bigskip

\vspace{15 pt} {H\"{U}SEYIN AYDIN,  Atat\"{u}rk University,
Faculty of Art and Sciences, Department of Mathematics, Erzurum,
Turkey. E-mail: aydinh@atauni.edu.tr}

\vspace{15 pt} {INCI G\"{U}LTEKYN,  Atat\"{u}rk University,
Faculty of Art and Sciences, Department of Mathematics, Erzurum,
Turkey. E-mail: inciakarg@yahoo.com}

\vspace{15 pt} {MICHELE MULAZZANI, Department of Mathematics and
C.I.R.A.M., University of Bologna, Bologna, Italy. E-mail:
mulazza@dm.unibo.it}

\end{document}